\documentclass[12pt,reqno]{amsart}
\usepackage{amssymb,delarray}
\usepackage{amsfonts}
\usepackage{epsfig}
\usepackage[all]{xy}

\def\leq{\leqslant}
\def\geq{\geqslant}

\newtheorem{thm}{Theorem}
\newtheorem{lem}
{Lemma}
\newtheorem{prop}
{Proposition}
\newtheorem{claim}
{Claim}
\newtheorem{df}[thm]{Definition}
\newtheorem{cor}
{Corollary}
\newtheorem{rem}
{Remark}
{Question}

{\catcode`\@=11
\gdef\n@te#1#2{\leavevmode\vadjust{%
 {\setbox\z@\hbox to\z@{\strut#1}%
  \setbox\z@\hbox{\raise\dp\strutbox\box\z@}\ht\z@=\z@\dp\z@=\z@%
  #2\box\z@}}}
\gdef\leftnote#1{\n@te{\hss#1\quad}{}}
\gdef\rightnote#1{\n@te{\quad\kern-\leftskip#1\hss}{\moveright\hsize}}
\gdef\?{\FN@\qumark}
\gdef\qumark{\ifx\next"\DN@"##1"{\leftnote{\rm##1}}\else
 \DN@{\leftnote{\rm??}}\fi{\rm??}\next@}}

\begin{document}
\baselineskip=13.7pt plus 2pt 

\title[The ambiguity index of an equipped finite group]
{The ambiguity index of an equipped finite group}
\author[F.A. Bogomolov and Vik.S. Kulikov ]
{F.A. Bogomolov and Vik.S. Kulikov }

\address{Courant Institute of Mathematical Sciences and
National Research University Higher School of Economics}
\email{bogomolo@cims.nyu.edu}

\address{Steklov Mathematical Institute}
\email{kulikov@mi.ras.ru}

\dedicatory{} \subjclass{}
\thanks{ The second author was partially supported
by grants of NSh-2998.2014.1, RFBR 14-01-00160, and both authors
were supported by AG Laboratory HSE, RF government grant, ag.
11.G34.31.0023. }

\begin{abstract}
In \cite{Ku0}, the ambiguity index $a_{(G,O)}$ was introduced for
each equipped finite group $(G,O)$. It is equal
to the number of connected components of a Hurwitz space
parametrizing
coverings of a projective line with Galois group $G$ assuming that
all local monodromies belong to conjugacy classes $O$ in $G$
and the number of branch points is greater than some constant.
 We prove in this article that
the ambiguity index can be identified with the size of a
generalization of so called Bogomolov multiplier (\cite{Kun1}, see also
\cite{BO87}) and hence can be easily computed for many pairs $(G,O)$.

\end{abstract}

\maketitle
\setcounter{tocdepth}{1}


\def\st{{\sf st}}



\section*{Introduction} \label{introduc}
Let $G$ be a finite group and $O$ be a subset of $G$ consisting of
conjugacy classes $C_i$ of $G$, $O= C_1\cup\dots\cup C_m$, which
together generate $G$. The pair $(G,O)$ is called an {\it equipped
group} and $O$ is called an {\it equipment} of $G$. We fix the
numbering of conjugacy classes  contained in $O$. One can associate
a  $C$-group $(\widetilde G, \widetilde O)$ to each equipped group
$(G,O)$. The $C$-group $\widetilde G$ is generated by the letters of
the alphabet $Y=Y_O=\{ y_g \mid g\in O\}$ subject to relations:

\begin{equation} \label{rel2} y_{g_1}y_{g_2}=
y_{g_2}y_{g_2^{-1}g_1g_2}=y_{g_1g_2g_1^{-1}} y_{g_1}.\end{equation}

We assume $\widetilde O=Y_O$ in the definition of $\widetilde G$  .

There is an obvious natural homomorphism $\beta :\widetilde G\to G$
given by  $\beta (y_g)=g$. It was shown in  \cite{Ku0},  that the
commutator subgroup $[\widetilde G,\widetilde G]$ is finite. The
order $a_{(G,O)}$ of the group $\ker \beta \cap [\widetilde
G,\widetilde G]$ was called the {\it ambiguity index} of the
equipped finite group $(G,O)$.

The notion of equipped groups is related to the description of
Hurwitz spaces parametrizing maps between projective curves with $G$
as  the monodromy group and the ambiguity index $a_{G,O}$ is equal
to the properly defined "asymptotic" number of connected components
of Hurwitz space parametrizing covering of curves with fixed
ramification data. More precisely, let $f:E\to F$ be a morphism
of a non-singular complex irreducible projective curve $E$ onto
a non-singular projective curve $F$. Let
us choose a point $z_0\in F$ such that $z_0$ is not a branch point
of $f$ hence the points  $f^{-1}(z_0)=\{ w_1,\dots ,w_d\}$, where
$d=\deg f$, are simple. If we fix the numbering of points in
$f^{-1}(z_0)$  then we call $f$ a  {\it marked covering}.

 Let $B=\{ z_1, \dots, z_n\}\subset F$ be the set of branch points of
$f$. The numbering of the points of  $f^{-1}(z_0)$ defines a
homomorphism $f_*:\pi_1(F\setminus B, z_0)\to \Sigma_d$ of the
fundamental group $\pi_1=\pi_1(F\setminus B, z_0)$ to the symmetric
group $\Sigma_d$. Define $G\subset \Sigma_d$ as $\text{im} f_* =G$.
It acts transitively on $f^{-1}(z_0)$.  Let
$\gamma_1,\dots,\gamma_n$ be simple loops around, respectively, the
points $z_1,\dots, z_n$ starting at $z_0$. The image
$g_j=f_*(\gamma_j)\in G$ is called a {\it local monodromy} of $f$ at
the point $z_j$. Each local monodromy $g_j$ depends on the choice of
$\gamma_j$, therefore it is defined uniquely up to conjugation in
$G$.

Denote by $O=C_1\cup\dots \cup C_m\subset G$ the union of conjugacy
classes of all local monodromies and by $\tau_i$ the number of local
monodromies of $f$ belonging to the conjugacy class $C_i$. The
collection $\tau=(\tau_1C_1,\dots,\tau_mC_m)$ is called the {\it
monodromy type} of $f$. Assume that the elements of $O$ generate the
group $G$. Then the pair $(G,O)$ is an equipped group.

Let $\text{HUR}_{d,G,O,\tau}^m(F,z_0)$ be the Hurwitz space (see the
definition of Hurwitz spaces in \cite{F} or in \cite{Ku-Kh}) of
marked degree $d$ coverings of $F$ with Galois group $G\subset
\Sigma_d$, local monodromies in $O$, and monodromy type $\tau$.
Hurwitz space $\text{HUR}_{d,G,O,\tau}^m(F,z_0)$ may consists of a
different number of connected components. However it was proved in
\cite{Ku-Kh} that for each equipped finite group $(G,O)$,
$O=C_1\cup\dots \cup C_m$, there is a number $T$ such that the
number of irreducible components of each non-empty Hurwitz space
$\text{HUR}^m_{d,G,O, \tau}(F, z_0)$ is equal to $a_{(G,O)}$ if
$\tau_i\geq T$ for all $i=1,\dots,m$. The number $T$ does not depend
on the base curve $F$ and degree $d$ of the coverings.

The subgroup $B_0(G)\subset H^2(G,Q/Z)$ was defined and studied in
\cite{BO87}. It consists of elements of  $H^2(G,Q/Z)$ which restrict
trivially onto abelian subgroups of $G$. It was conjectured in
\cite{BMP} that $B_0(G)$ is trivial for simple groups. This
conjecture was partially solved  already in \cite{BMP}
and it was  completely solved by Kunyavski in  \cite{Kun1},
and by Kunyavski-Kang in  \cite{Kun} for a wider class of almost simple groups. The
latter consists of groups $G$ which contain some simple group $L$
and in turn are contained in the automorphism group $Aut L$.
Kunyavski in \cite{Kun1} called $B_0(G)$ as {\it Bogomolov
multiplier} and we are going to use his terminology here. Denote by
$b_0(G)$ the order of the group $B_0(G)$ and denote by $h_2(G)$ the
Schur multiplier of the group $G$, that is, the order of the group
$H_2(G,\mathbb Z)$.

The aim of this article is to prove
\begin{thm} \label{main} For an  equipped finite group $(G,O)$ we have
the following inequalities $$b_0(G)\leq a_{(G,O)}\leq h_2(G).$$ In
particular, $a_{(G,G\setminus \{1\})}=b_0(G)$.
\end{thm}

Since, by \cite{Kun1}, $b_0(G)=1$ for a finite almost simple group
$G$, we conclude:

\begin{cor} \label{cor1} Let $G$ be a finite almost simple group.
Then there is a constant $T$ such that for any projective
irreducible non-singular curve $F$ each non-empty Hurwitz space
$\text{HUR}^m_{d,G,G\setminus \{ 1\}, \tau}(F, z_0)$ is irreducible
if all $\tau_i\geq T$.
\end{cor}

It was shown in \cite{Ku0} that if $O_1\subset O_2$ are two
equipments of a finite group $G$, then $a_{(G,O_2)}\leq
a_{(G,O_1)}$.

For a symmetric group $\Sigma_d$, the famous
Clebsch -- Hurwitz Theorem (\cite{Cl}, \cite{H}) states that the
ambiguity index $a_{(\Sigma_d,T)}=1$, where $T$ is the set of
transpositions in $\Sigma_d$, and it was shown in \cite{Ku2} that
the ambiguity index $a_{(\Sigma_d,O)}=1$ if the equipment $O$
contains an odd permutation $\sigma\in \Sigma_d$ such that $\sigma$
leaves fixed at least two elements.   Theorem \ref{sym} (see
subsection \ref{subsym}) gives the complete answer on the value of
$a_{(\Sigma_d,O)}$ for each equipment $O$ of $\Sigma_d$.  Also in
subsection \ref{subsym}, we give the complete answer on the value
of $a_{(\mathbb A_d,O)}$ for each $d$ and for each equipment $O$ of
the alternating group $\mathbb A_d$.

In Section \ref{sec1}, we remind some properties of $C$-groups and
prove one of the inequalities claimed in Theorem \ref{main}. In
Section \ref{sec2}, we complete the proof of this Theorem.

In Section \ref{sec3}, we investigate the properties of ambiguity
indices of a quasi-cover of an equipped finite group $(G,O)$ and in
Section \ref{sec4}, we give a cohomologycal description of the
ambiguity indices.

In Section \ref{sec5}, we give examples of finite groups $G$ which
Bogomolov multiplier $b_0(G)>1$. Therefore for such groups $G$ each
non-empty 
space $\text{HUR}^m_{d,G,O, \tau}(F, z_0)$
consists of at least $b_0(G)>1$ irreducible components for any
$\tau=(\tau_1,\dots,\tau_m)$ with big enough $\tau_i$.

In this article, if $\mathbb F$ is a free group freely generated by
an alphabet $X$, $N$ is a normal subgroup of $\mathbb F$, and a
group $G=\mathbb F/N$, then a word $w=w(x_{i_1},\dots, x_{i_n})$ in
letters $x_{i_j}\in X$ and their inverses will be considered as an
element of $G$ in case if it does not lead to misunderstanding.

\section{$C$-groups and their properties} \label{sec1}
Let us remind the definition of a $C$-group (see, for example,
\cite{Ku}).

\begin{df} A group $G$ is a $C$-group if there is a set of generators
$x\in X$ in $G$ such that the basis of relations between $x\in X$
consists of the following relations:
\begin{equation} \label{C-rel} x_i^{-1}x_jx_i=x_k,\quad (x_i,x_j,x_k)\in M,
\end{equation}
where $M$ is a subset of $X^3$.
\end{df}

Thus the $C$-structure of $G$ is defined by $X\subset G$ and
$M\subset X^3$.

Let $\mathbb F$ be a free group freely generated by an alphabet $X$.
Denote by $N$ the subgroup of $\mathbb F$ normally generated by the
elements $x_i^{-1}x_jx_ix^{-1}_k$, $(x_i,x_j,x_k)\in M$. The group
$N$ is a normal subgroup of $\mathbb F$. Let $f: \mathbb F\to
G=\mathbb F/N$ be the natural epimorphism given by presentation
(\ref{C-rel}).  In the sequel, we consider each $C$-group $G$ as an
equipped group $(G,O)$ with the equipment $O=f(X^{\mathbb F})$
(where $X^{\mathbb F}$ is the orbit of $X$ under the action of the
group of inner automorphisms of $\mathbb F$). The elements of $O$
are called  $C$-{\it generators} of the $C$-group $G$. In
particular, the equipped group $(\mathbb F,X^{\mathbb F})$ is a
$C$-group.

A homomorphism $f:G_1\to G_2$ of a $C$-group $(G_1,O_1)$ to a
$C$-group $(G_2,O_2)$ is called a $C$-{\it homomorphism} if it is a
homomorphism of equipped groups, that is, $f(O_1)\subset O_2$. In
particular, two $C$-groups $(G_1,O_1)$ and $(G_2,O_2)$ are $C$-{\it
isomorphic} if they are isomorphic as equipped groups.

\begin{claim} ({\rm Lemma 3.6 in} \cite{Ku}) \label{cl1} Let $N$ be
a normal subgroup of $\mathbb F$ normally generated by elements of
the form $w_i^{-1}x_jw_iw_lx_k^{-1}w_l^{-1}$, where $w_i$ and $w_l$
are some elements of $\mathbb F$ and $x_j,x_k\in X$. Let $f:\mathbb
F\to G\simeq \mathbb F/N$ be the natural epimorphism. Then
$(G,f(X^{\mathbb F}))$ is a $C$-group and $f$ is a $C$-homomorphism.
\end{claim}

To each $C$-group $(G,O)$, one can associate a $C$-graph. By
definition, the $C$-{\it graph $\Gamma=\Gamma_{(G,O)}$ of a
$C$-group} $(G,O)$ is a directed  labeled graph whose set of
vertices $V=\{ v_{g_i}\mid g_i\in O\}$ is in one to one
correspondence with the set $O$. Two vertices $v_{g_1}$ and
$v_{g_2}$, $g_1, g_2\in O$, are connected by a labeled edge
$e_{v_{g_1},v_{g_2},v_g}$ (here $v_{g_1}$ is the tail of
$e_{v_{g_1},v_{g_2},v_g}$, $v_{g_2}$ is the head of
$e_{v_{g_1},v_{g_2},v_g}$, and $v_{g}$ is the label of
$e_{v_{g_1},v_{g_2},v_g}$) if and only if in $G$ we have the
relation $g^{-1}g_1g=g_2$ with some $g\in O$.

A $C$-homomorphism $f:(G_1,O_1)\to (G_2,O_2)$ of $C$-groups induces
a  map $f_*:\Gamma_{(G_1,O_1)}\to \Gamma_{(G_2,O_2)}$ from the
$C$-graph $\Gamma_{(G_1,O_1)}$ in the $C$-graph
$\Gamma_{(G_2,O_2)}$, where by definition, $f_*(v_g)=v_{f(g)}$ for
each vertex $v_g$ of $\Gamma_{(G_1,O_1)}$ and
$f_*(e_{v_{g_1},v_{g_2},v_g})=e_{v_{f(g_1)},v_{f(g_2)},v_{f(g)}}$
for each edge $e_{v_{g_1},v_{g_2},v_g}$ of $\Gamma_{(G_1,O_1)}$.

The following Claim is obvious.
\begin{claim} \label{cl2} A $C$-homomorphism
$f:(G_1,O_1)\to (G_2,O_2)$ is a $C$-isomorphism if $f_*$ is
one-to-one between the sets of vertices of $\Gamma_{(G_1,O_1)}$ and
$\Gamma_{(G_2,O_2)}$.
\end{claim}

In the sequel, we will consider only finitely presented $C$-groups
(as groups without equipment) and $C$-graphs consisting of finitely
many connected components. Denote by $m$ the number of connected
components of  a $C$-graph $\Gamma_{(G,O)}$.

Then it is easy to see that $G/[G,G]\simeq \mathbb Z^m$ and any two
$C$-generators $g_1$ and $g_2$ are conjugated in the $C$-group $G$
if and only if $v_{g_1}$ and $v_{g_2}$ belong to the same connected
component of $\Gamma_{(G,O)}$.

Thus the set $O$ of $C$-generators of the $C$-group $(G,O)$ is the
union of $m$ conjugacy classes of $G$ and there is a one-to-one
correspondence between the conjugacy classes of $G$ contained in $O$
and the set of connected components of $\Gamma_{(G,O)}$.

Denote by $\tau: G_{\Gamma}\to H_1(G,\mathbb Z)=G/[G,G]$ the natural
epimorphism. In the sequel, we fix some numbering of the connected
components of $\Gamma_{(G,O)}$. Then the group $H_1(G,\mathbb
Z)\simeq \mathbb Z^m$ obtains a natural base consisting of vectors
$\tau(g)=(0,\dots,0,1,0\dots,0)$, where $1$ stands on the $i$-th
place if $g$ is a $C$-generator of $G$ and $v_g$ belongs to the
$i$-th connected component of $\Gamma_{(G,O)}$. For $g\in G$ the
image $\tau(g)$ is called the {\it type} of $g$.

\begin{lem} \label{lem2} Let $g_1, g_2$ be two $C$-generators
of a $C$-group $(G,O)$, $N$ the normal closure of $g_1g_2^{-1}$ in
$G$, and $f: G\to G_1=G/N$ the natural epimorphism. Then
\begin{itemize}
\item[($i$)] $(G_1,O_1)$ is a $C$-group, where $O_1=f(O)$,
and $f$ is a $C$-homomorphism;
\item[($ii$)] the map $f_*:\Gamma_{(G,O)}\to \Gamma_{(G_1,O_1)}$ is
a surjection.
\item[($iii$)] if
$g_1g_2^{-1}$ belong to the center $Z(G)$ of the group $G$ and let
$v_{g_1}$ and $v_{g_2}$ belong to different components of
$\Gamma_{(G,O)}$, then
\begin{itemize}
\item[($iii_1$)] the number of connected components of the $C$-graph
$\Gamma_{(G_1,O_1)}$ is less than the number of connected components
of the $C$-graph $\Gamma_{(G,O)}$;
\item[($iii_2$)] $f: [G,G]\to [G_1,G_1]$ is an isomorphism.
\end{itemize}
\end{itemize}
\end{lem}
\proof Claims ($i$), ($ii$), and $(iii_1$) are obvious.

To prove $(iii_2)$, note that $N$ is a cyclic group generated by
$g_1g_2^{-1}$, since $g_1g_2^{-1}$ belongs to the center $Z(G)$. The
type $\tau((g_1g_2^{-1})^n)$ is non-zero for $n\neq 0$, since
$v_{g_1}$ and $v_{g_2}$ belong to different connected components of
$\Gamma_{(G,O)}$. Therefore to complete the proof, it suffices to
note that the groups $N$ and $[G,G]$ have trivial intersection,
since $\tau (g)=0$ for all $g\in [G,G]$. \qed \\

A $C$-group $(G,O)$ is called a  $C$-{\it finite group} if the set
of vertices of $C$-graph $\Gamma_{(G,O)}$ is  finite or, the same,
if the equipment $O$ of $G$ is a finite set.

\begin{prop} {\rm (\cite{Ku0})} \label{fincom} Let $(G,O)$ be
a $C$-finite group. Then the commutator $[G,G]$ is a finite group.
\end{prop}

As it was mention in Introduction, for each finite equipped group
$(G,O)$, one can associate a $C$-group $(\widetilde G, \widetilde
O)$ defined as follows. The group $\widetilde G$ is generated by the
letters of the alphabet $Y=Y_O=\{ y_g \mid g\in O\}$ subject to
relations
\begin{equation} \label{rel3} y_{g_1}y_{g_2}
=y_{g_2}y_{g_2^{-1}g_1g_2}=y_{g_1g_2g_1^{-1}} y_{g_1}.\end{equation}
Here $\widetilde O=Y_O$ and  there is a natural epimorphism
$\beta_O :\widetilde G\to G$ given by $\beta_O (y_g)=g$.

Note also that a homomorphism of equipped groups $f:(G_1,O_1)\to
(G,O)$ induces a $C$-homomorphism $\widetilde f:(\widetilde
G_1,\widetilde O_1)\to (\widetilde G,\widetilde O)$ such that
$f\circ \beta_{O_1}=\widetilde f\circ \beta_O$.

Let the elements of a subset $S$ of an equipment $O$ of a group $G$
generate the group $G$ and $O=S^G$. Denote by $\mathbb F_S$ a free
group freely generated by the alphabet $Y_S=\{ y_g \mid g\in S\}$
and $R_S$ is the normal subgroup of $\mathbb F_S$ such that the
natural epimorphism $h_S: \mathbb F_S\to \mathbb F_S/R_S\simeq G$
gives a presentation of the group $G$.

\begin{claim} \label{cl3} Let $\widetilde R_S\subset R_S$ be the
normal subgroup normally generated by the elements of $R_S$ of the
form $w_{i,j}^{-1}y_{g_i}w_{i,j}y_{g_j}^{-1}$, where $w_{i,j}\in
\mathbb F_S$ and $y_{g_i},y_{g_j}\in Y_S$. Then the $C$-group
$(\widetilde G,\widetilde O)$ has the following presentation:
$\widetilde G\simeq \mathbb F_S/\widetilde R_S$ such that the images
of the elements of $Y_S$ are $C$-generators of $\widetilde G$.
\end{claim}

\proof Denote by $G_1=\mathbb F_S/\widetilde R_S$. By Claim
\ref{cl1},  $G_1$ is a $C$-group with $C$-equipment $O_1=Y_S^{G_1}$
and there is a natural epimorphism $\beta_S :((G_1,O_1)\to (G,O)$
given by $\beta_S(y_g)=g$ for $g\in S$.

Assume that $S$ consists of elements $g_1,\dots,g_n\in O$. If $S\neq
O$ then  choose an element $g_{n+1}\in O\setminus S$. It is
conjugated to some $g_i\in S$. Denote by $R_{g_{n+1}}$ the set of
all presentations of $g_{n+1}$ in the form
\begin{equation} \label{rele1} g_{n+1}=w(g_1,\dots,g_n)^{-1}gw(g_1,\dots,g_n),
\quad g\in S.\end{equation} Note that if
$$g_{n+1}=w_i(g_1,\dots,g_n)^{-1}g_iw_i(g_1,\dots,g_n)\, \,
\text{and}\, \,
g_{n+1}=w_j(g_1,\dots,g_n)^{-1}g_jw_j(g_1,\dots,g_n),$$ then
$$w_jw_i^{-1}g_iw_iw_j^{-1}=g_j,$$ that is,
\begin{equation} \label{rele2} w_j(y_{g_1},\dots,y_{g_n})w_i^{-1}(y_{g_1},
\dots,y_{g_n})y_{g_i}w_i(y_{g_1},\dots,y_{g_n})
w_j^{-1}(y_{g_1},\dots,y_{g_n})y_{g_j}^{-1}\in R_S.\end{equation}
Similarly, if $g_{n+1}=w_i(g_1,\dots,g_n)$ and
$g_{n+1}^{-1}g_ig_{n+1}=g_j$ for some $g_i,g_j\in S$, then
\begin{equation} \label{rele3}
w(y_{g_1},\dots,y_{g_n})^{-1}y_{g_i}w(y_{g_1},\dots,y_{g_n})y_{g_j}^{-1}\in
R_S.\end{equation}

Therefore, if $S_1=S\cup \{ g_{n+1}\}$, $\mathbb F_{S_1}$ is a free
group freely generated by the alphabet $Y_{S_1}=\{ y_g  \mid g\in
S_1\}$, $R_{g_{n+1}}$ is the set of words of the form
$$w(y_{g_1},\dots,y{g_n})^{-1}y_gw(y_{g_1},\dots,y_{g_n})y_{g_{n+1}}^{-1}$$
defined by all relations (\ref{rele1}), and $\widetilde R_{S_1}$ is
the normal closure in $\mathbb F_{S_1}$ of the set $\widetilde
R_S\cup R{g_{n+1}}$, then $G_1\simeq \mathbb F_{S_1}/\widetilde
R_{S_1}$ in view of relations (\ref{rele2}) and (\ref{rele3}).

Note that if we have a relation $g_i^{-1}g_jg_i=g_k$ for some
$g_i,g_j,g_k\in S_1$ then
\begin{equation} \label{rele4} y_{g_i}^{-1}y_{g_j}y_{g_i}y_{g_k}^{-1}\in
 \widetilde R_{S_1}. \end{equation}

If $S_1\neq O$, then we can repeat the construction described above
and obtain a presentation $G_1\simeq \mathbb F_{S_2}/\widetilde
R_{S_2}$, and so on. After several steps we obtain a presentation
$G_1\simeq \mathbb F_{O}/\widetilde R_{O}$. Note that, by induction,
we deduce that for any relation in $G$ of the form
$g_i^{-1}g_jg_i=g_k$ for some $g_i,g_j,g_k\in O$ we have
$y_{g_i}^{-1}y_{g_j}y_{g_i}y_{g_k}^{-1}\in \widetilde R_{O}$.
Therefore there is a natural $C$-homomorphism $f:(\widetilde
G,\widetilde O)\to (G_1,O_1)$.
By Claim \ref{cl2}, $f$ is a $C$-isomorphism. \qed \\

For an equipped finite group $(C,O)$, consider a presentation of $G$
of the following form. Let us take a free group $\mathbb F=\mathbb
F_O$ freely generated by the alphabet $X_O=\{ x_g \mid g\in O\}$.
Consider a normal subgroup $R_O\subset \mathbb F $ such that
$\mathbb F/R_O\simeq G$. Let $h_O: \mathbb F\to \mathbb F/R_O\simeq
G$ be the natural epimorphism.

We can associate to $(G,O)$ a group $\overline G =\mathbb F/[\mathbb
F,R_O]$.

Denote by  $\alpha_O :\overline G \to G$ the natural epimorphism. By
Claim \ref{cl1}, $(\overline G, \overline O)$ is a $C$- group, where
$\overline O= h_O(X_O^{\mathbb F})$. It is evident that there is the
natural epimorphism of $C$-groups $\kappa_O: (\overline G,\overline
O)\to (\widetilde G,\widetilde O)$ sending $\kappa_O(x_g)=y_g$ for
all $g\in O$ and such that $\alpha_O= \beta_O\circ \kappa_O$. The
$C$-group $(\overline G,\overline O)$ is called the {\it universal
central $C$-extension} of the equipped finite group $(G,O)$.

It is easy to see that $\alpha_O :\overline G\to G$ is a central
extension of groups, that is, $\ker \alpha_O$ is a subgroup of the
center $Z(\overline G)$.

We have $$\ker \alpha_O\cap [\overline G,\overline G]=(R_O\cap
[\mathbb F,\mathbb F])/ [\mathbb F,R_O].$$ By Hopf's integral
homology formula, we have $H_2(G,\mathbb Z)\simeq (R_O\cap[\mathbb
F,\mathbb F])/[\mathbb F,R_O]$. Denote by $h_2(G)$ the order of the
group $H_2(G,\mathbb Z)$ and denote by $K_{(G,O)}$ the subgroup of
$(R_O\cap[\mathbb F,\mathbb F])/[\mathbb F,R_O]$ generated by the
elements of $R_O$ of the form $[w,x_g]$, where $g\in O$ and $w\in
\mathbb F$, and let $k_{(G,O)}$ be its order.

\begin{thm} \label{thm2} For an equipped finite group $(G,O)$ we have
$$ h_2(G)=k_{(G,O)}a_{(G,O)}.$$
\end{thm}
\proof We have $\ker \kappa_O\subset \ker\alpha_O$. Therefore $\ker
\kappa_O\subset Z(\overline G)$.

Let us show that for some $n\geq 0$ there is a sequence of
$C$-groups $\overline G_0=\mathbb F/R_0,\dots ,\overline G_n=\mathbb
F/R_n$, a sequence of $C$-homomorphisms
$$\varphi_i:(\overline G_i,\overline O_i)\to (\overline  G_{i+1},
\overline O_{i+1}),\qquad 0\leq i\leq n-1,$$ where $(\overline
G_0,\overline O_0)=(\overline G,\overline O)$, and a
$C$-homomorphism $\overline{\kappa}: (\overline G_n,\overline
O_n)\to (\widetilde G,\widetilde O)$ such that

\begin{itemize}
\item[$(i)$] $\kappa=\overline{\kappa}\circ \varphi$,
where $\varphi=\varphi_{n-1}\circ\dots\circ \varphi_0$;
\item[$(ii)$]  for each $i$ the homomorphism
$\varphi_i:[\overline G_i,\overline G_i]\to [\overline
G_{i+1},\overline G_{i+1}]$ is an isomorphism;
\item[$(iii)$] $\overline{\kappa}_{*}$ induces a one-to-one
correspondence between the connected components of the $C$-graphs
$\Gamma_{(\overline G_n,\overline O_n)}$ and $\Gamma_{(\widetilde
G,\widetilde O)}$.
\end{itemize}

Indeed, let us put $R_0=R_O$ and consider the map $\kappa_*$. If it
is induces a one-to-one correspondence between the connected
components of the $C$-graphs $\Gamma_{(\overline G,\overline O)}$
and $\Gamma_{(\widetilde G,\widetilde O)}$, then $n=0$ and it is
nothing to prove.

Otherwise, for some $g\in O$ there is a vertex $v_{y_g}$ of
$\Gamma_{(\widetilde G,\widetilde O)}$ which preimage
$\kappa_*^{-1}(v_{y_g})$ contains at least two vertices, say
$v_{x_g}$ and $v_{\overline g}$ (here $\overline g$ is an element of
$X^{\mathbb F}$), of $\Gamma_{(\overline G,\overline O)}$ belonging
to different connected components of $\Gamma_{(\overline G,\overline
O)}$.

Denote by $R_1$ the normal closure of $R_O\cup \{ x_g\overline
g^{-1}\}$ in $\mathbb F$  and consider the natural homomorphism
$\varphi_0:\overline G\to \overline G_1=\mathbb F/R_1$. The element
$x_g\overline g^{-1}$, considered as an element of $\overline G$,
belongs to $\ker \kappa$. Therefore, $x_g\overline g^{-1}\in
Z(\overline G)$.

Denote by $\kappa_1:\overline G_1\to\widetilde G$ the homomorphism
induced by $\kappa$. By Lemma \ref{lem2}, the homomorphism
$\varphi_1$ is a $C$-homomorphism of $C$-groups. It is easy to see
that $\varphi_0: [\overline G_0,\overline G_0]\to [\overline
G_{1},\overline G_{1}]$ is an isomorphism and the number of
connected components of the $C$-graph $\Gamma_{(\overline
G_1,\overline O_1)}$ is less than the number of connected components
of the $C$-graph $\Gamma_{(\overline G,\overline O)}$.

Assume now that $\kappa_{1*}$ is not a one-to-one correspondence
between the connected components of the $C$-graphs
$\Gamma_{(\overline G1,\overline O1)}$ and $\Gamma_{(\widetilde
G,\widetilde O)}$. Then for some $g_1\in O$ there is a vertex
$v_{y_{g_1}}$ of $\Gamma_{(\widetilde G,\widetilde O)}$ which
preimage $\kappa_{1*}^{-1}(v_{y_{g_1}})$ contains at least two
vertices $v_{x_{g_1}}$ and $v_{\overline g_1}$ of
$\Gamma_{(\overline G_1,\overline O_1)}$ belonging to different
connected components of $\Gamma_{(\overline G_1,\overline O_1)}$.

Hence we can repeat the construction described above and obtain a
$C$-group $(\overline G_2,\overline O_2)$ and $C$-homomorphisms
$\varphi_1:\overline G_1\to \overline G_2=\mathbb F/R_2$ and
$\kappa_2:\overline G_2\to\widetilde G$ such that $\varphi_1:
[\overline G_1,\overline G_1]\to [\overline G_{2},\overline G_{2}]$
is an isomorphism and the number of connected components of the
$C$-graph $\Gamma_{(\overline G_2,\overline O_2)}$ is less than the
number of connected components of the $C$-graph $\Gamma_{(\overline
G_1,\overline O_1)}$. Since the number of connected components of
the $C$-graph $\Gamma_{(\overline G,\overline O)}$ is finite, after
several ($n$) steps of our construction we obtain the desired
sequences of $C$-groups and $C$-homomorphisms.

Now, consider the $C$-homomorphism $\overline{\kappa} : \overline
G_n\to \widetilde G$. The $C$-graph $\Gamma_{(\widetilde
G,\widetilde O)}$ consists of connected components $\Gamma_1,\dots
,\Gamma_m$ . Let $\{v_{g_{i,1}},\dots, v_{g_{i,l_i}}\}$ be the set
of the vertices of $\Gamma_i$. We have $O=\{g_{i,j}\}_{1\leq i\leq
m, 1\leq j \leq l_i}$. Then
$\overline{\Gamma}_i=\overline{\kappa}_*^{-1}(\Gamma_i)$ are the
connected components of $\Gamma_{(\overline G_n,\overline O_n)}$.

Let $\overline{\kappa}_n^{-1}(v_{y_{g_{i,j}}})=\{v_{x_{g_{i,j}}},
v_{\overline g_{i,j,1}},\dots,  v_{\overline g_{i,j,r_{i,j}}}\}$,
$\overline g_{i,j,k}\in \overline O_n$ for $1\leq k\leq r_{i,j}$.

Since the graph $\overline{\Gamma}_i$ is connected, there are words
$w_{i,j,k}$ in letters of $X_O$ and their inverses such that
$$\overline g_{i,j,k}=w_{i,j,k}x_{g_{i,j}}w_{i,j,k}^{-1},\, \, \, \,
\, \,  1\leq k\leq r_{i,j}.$$ Obviously, the elements
$u_{i,j,k}=[w_{i,j,k},x_{g_{i,j}}]=\overline
g_{i,j,k}x_{g_{i,j}}^{-1}$ belong to $[\overline G_n,\overline
G_n]\cap\ker \overline{\kappa}$. Therefore $u_{i,j,k}$, as elements
of $\mathbb F$ belong to $R_O\cap [\mathbb F,\mathbb F]$.

Consider the group $\overline G_{n+1}=\mathbb F/R_{n+1}$, where the
group $R_{n+1}$ is the normal closure of $R_n\cup
\{u_{i,j,k}\}_{1\leq i\leq m, 1\leq j\leq l_i, 1\leq k\leq r_{i,j}}$
in $\mathbb F$. Then, by Claim \ref{cl1}, $\overline G_{n+1}=\mathbb
F/R_{n+1}$ is a $C$-group and the natural map
$\overline{\kappa}_{1}: \overline G_{n+1}\to \widetilde G$, induced
by $\overline{\kappa}$, is a $C$-homomorphism. Moreover, $\ker
\varphi_{n}$ of the natural epimorphism $\varphi_{n}: \overline
G_n\to\overline G_{n+1}$  is a subgroup of $[\overline G_n,\overline
G_n]\simeq [\overline G,\overline G]=[\mathbb F,\mathbb F]/[\mathbb
F,R_O]$ generated by the elements
$u_{i,j,k}=[w_{i,j,k},x_{g_{i,j}}]$, where $1\leq i\leq m$, $1\leq
j\leq l_i$, and $1\leq k\leq r_{i,j}$.

To complete the proof of Theorem \ref{thm2}, it suffices to note
that $\overline{\kappa}_{1*}$ induces a one-to-one correspondence
between the sets of vertices of the $C$-graphs  $\Gamma_{(\overline
G_{n+1},\overline O_{n+1})}$ and $\Gamma_{(\widetilde G,\widetilde
O)}$, since all $u_{i,j,k}=\overline g_{i,j,k}x_{g_{i,j}}^{-1}$
belong to $\ker \varphi_n$. Therefore  $\overline{\kappa}_1$ is an
isomorphism. \qed

\begin{lem} \label{order} Let the order of $g\in O$
be $n$ and let  
$[x_g,w]\in ([\mathbb F,\mathbb F]\cap R_O)/[\mathbb F, R_O]
\subset\mathbb F/[\mathbb F, R_O]$. Then the order of the element
$[x_g,w]$ is a divisor of $n$.
\end{lem}
\proof  The elements $x_g^n$ and $[x_g,w]$ belong to the center of
the group $\mathbb F/[\mathbb F, R_O]$. Therefore
$$[x_g^n,w]=x_g^{n-1}[x_g,w]x_g^{1-n}[x_g^{n-1},w]=[x_g,w][x_g^{n-1},w]=\dots= [x_g,w]^n$$
is the unity of $\mathbb F/[\mathbb F, R_O]$. \qed

\begin{prop} \label{prop1} Let the equipment $O$ of an equipped finite group
$(G,O)$ consists of conjugacy classes of elements of orders coprime
with $h_2(G)$. Then $a_{(G,O)}=h_2(G)$.
\end{prop}
\proof  It follows from Lemma \ref{order} and Theorem \ref{thm2}.
\qed

\section{Proof of Theorem \ref{main}} \label{sec2}
By definition, the {\it Bogomolov multiplier} $b_0(G)$ of a finite
group $G$ is the order of the group $$B_0(G)=\ker[H^2(G,\mathbb
Q/\mathbb Z)\to \bigotimes_{A\subset G}H^2(A,\mathbb Q/\mathbb Z)]$$
where $A$ runs over all abelian subgroups of $G$.

\begin{rem} {\rm Note that it suffices to consider only restrictions to
abelian groups with two generators in order to define that the
element $w\in H^2(G,\mathbb Q/\mathbb Z) $ is contained in
$B_0(G)$.}
\end{rem}

There is a natural duality between $H^2(G,\mathbb Q/\mathbb Z)$ and
$H_2(G,\mathbb Z)$. Both groups are finite for finite groups $G$ and
duality implies an isomorphism of $H^2(G,\mathbb Q/\mathbb Z)$ and
$Hom (H_2(G,\mathbb Z),\mathbb Q/\mathbb Z)$ as abstract groups.

By Theorem \ref{thm2}, we have the inequality $h_2(G)\geq a_{(G,O)}$
for any equipped finite group $(G,O)$. By Corollary 2 in \cite{Ku0},
we have inequality $a_{(G,O)}\geq a_{(G,G\setminus \{ 1\})}$ for
each equipment $O$ of $G$. Therefore to prove Theorem \ref{main} it
suffice to show that for the equipped finite group $(G,G\setminus \{
1\})$ its ambiguity index $a_{(G,G\setminus \{ 1\})}$ is equal to
$b_0(G)$.

In notation used in Section \ref{sec1} and by Theorem \ref{thm2}, we
have $$a_{(G,G\setminus \{ 1\})}=\frac{h_2(G)}{k_{(G,G\setminus \{
1\})}},$$ where $k_{(G,G\setminus \{ 1\})}$ is the order of the
subgroup $K_{G\setminus \{ 1\}}$ of the group $$(R_{G\setminus \{
1\}}\cap[\mathbb F_{G\setminus \{ 1\}},\mathbb F_{G\setminus \{
1\}}])/[\mathbb F_{G\setminus \{ 1\}},R_{G\setminus \{ 1\}}]\simeq
H_2(G,\mathbb Z)$$ generated by the elements of $R_{G\setminus \{
1\}}$ of the form $[w,x_g]$, where $g\in G\setminus \{ 1\}$ and
$w\in \mathbb F_{G\setminus \{ 1\}}$.

\begin{lem} \label{generations} Let for some $w_1,w_2\in \mathbb F_{G\setminus \{
1\}}$ the commutator $[w_1,w_2]$ belong to $R_{G\setminus \{ 1\}}$.
Then $[w_1,w_2]$, considered as an element of $\mathbb F_{G\setminus
\{ 1\}}/[\mathbb F_{G\setminus \{ 1\}},R_{G\setminus \{ 1\}}]$,
belongs to $K_{G\setminus \{ 1\}}$.
\end{lem}

\proof First of all, note that if $[x_g, w]\in K_{G\setminus \{
1\}}$, then $[x_g, w]=[w,x_g^{-1}]=[x_g^{-1},w^{-1}]=[x_g^{-1}, w]$
in $K_{G\setminus \{ 1\}}$, since $K_{G\setminus \{ 1\}}$ is a
subgroup of the center of the $C$-group $\overline G_{G\setminus \{
1\}} =\mathbb F_{G\setminus \{ 1\}}/[\mathbb F_{G\setminus \{
1\}},R_{G\setminus \{ 1\}}]$ and these four commutators are
conjugated to each other in $\mathbb F_{G\setminus \{ 1\}}$.
Similarly, $[w,x_g]=[x_g,w^{-1}]=[w^{-1},x_g^{-1}]=[x_g^{-1},
w^{-1}]\in K_{G\setminus \{ 1\}}$, since $[w,x_g]$ is the inverse
element to the element $[x_g,w]$. Note also that for any $w_1$ the
element $w_1[w,x_{g}]w_1^{-1}$ belongs to $K_{G\setminus \{ 1\}}$ if
$[w,x_g]\in K_{G\setminus \{ 1\}}$.

Next, the elements $w_1^{-1}$ and $w_2^{-1}$, considered as elements
of $G$, are equal to some elements $g_1$ and $g_2$ of $G$. Therefore
if $[w_1,w_2]\in R_{G\setminus \{ 1\}}$ then
$$w_1x_{g_1},w_2x_{g_2},[x_{g_1},x_{g_2}],
[w_{2},x_{g_1}],[w_{1},x_{g_2}]\in R_{G\setminus \{ 1\}}.$$

In addition, we have $[w_1,w_2x_{g_2}]\in [\mathbb F_{G\setminus \{
1\}},R_{G\setminus \{ 1\}}]$ and
$$[w_1,w_2x_{g_2}]=[w_1,w_2](w_2[w_1,x_{g_2}]w_2^{-1}).$$
Therefore $[w_1,w_2]\in R_{G\setminus \{ 1\}}\cap[\mathbb
F_{G\setminus \{ 1\}},\mathbb F_{G\setminus \{ 1\}}]$ (as an element
of $K_{G\setminus \{ 1\}}$) is the inverse element to the element
$[w_1,x_{g_2}]\in K_{G\setminus \{ 1\}}$ and
hence $[w_1,w_{2}]\in K_{G\setminus \{ 1\}}$. \qed \\

To complete the proof of Theorem \ref{main}, note that, by Lemma
\ref{generations}, for each imbedding $i:H\to G$ of an abelian group
$H$ generated by two elements the image of $i_*:H_2(H,\mathbb Z)\to
H_2(G,\mathbb Z)$ is a subgroup of $K_{G\setminus \{ 1\}}$ and the
group $K_{G\setminus \{ 1\}}$ is generated by the images of such
elements. Therefore the group $$K^{\bot}_{G\setminus \{ 1\}}=\{
\varphi \in Hom (H_2(G,\mathbb Z),\mathbb Q/\mathbb Z)\mid
\varphi(w)=0\,\, \text{for all}\, \, w\in K_{G\setminus \{ 1\}}\}$$
coincides with the group $B_0(G)$ and its order is equal to
$a_{(G,G\setminus \{ 1\})}=\frac{h_2(G)}{k_{(G,G\setminus \{
1\})}}$. \qed

\section{Quasi-covers of equipped finite groups} \label{sec3}
In this section we use notations introduced in Sections \ref{sec1}.

\subsection{Definitions}
Let $f:(G_1,O_1)\to (G,O)$ be a homomorphism of equipped groups. We
say that $f$ is a {\it cover of equipped groups} (or, equivalently,
$(G_1,O_1)$ is a {\it cover} of $(G,O)$) if
\begin{itemize}
\item[($i$)] $f$ is an epimorphism such that $f(O_1)=O$;
\item[($ii$)] $\ker f$ is a subgroup of the center $ZG_1$ of $G_1$;
\item[($iii$)] $f_*: H_1(G_1,\mathbb Z)\to H_1(G,\mathbb Z)$ is
an isomorphism.
\end{itemize}

Let $f:(G_1,O_1)\to (G,O)$ be a homomorphism of equipped finite
groups. We say that $S\subset O_1$ is a {\it section} of $f$ if
$f_{|S}:S\to O$ is a one-to-one correspondence. Denote by
$O_S\subset O_1$ the orbit of $S$ under the action of the group of
the inner automorphisms of $G_1$.

Let $f:(G_1,O_1)\to (G,O)$ be an epimorphism of equipped groups such
that $\ker f\subset ZG_1$. We say that $f$ is a {\it quasi-cover of
equipped groups} (or, equivalently, $(G_1,O_1)$ is a {\it
quasi-cover} of $(G,O)$) if there is a section $S$ of $f$ such that
$O_S=O_1$.

Below, we will assume that for a quasi-cover $f$ of equipped groups
a section $S$ is chosen and fixed.

\subsection{Properties of quasi-covers}
\begin{lem} \label{section} Let $f:(G_1,O_1)\to (G,O)$ be a cover of
equipped finite groups and $S\subset O_1$ a section. Then $G_1$ is
generated by the elements of $S$.
\end{lem}

\proof Denote by $G_{S}$ the subgroup of $G_1$ generated by the
elements of $S$. Obviously, $\varphi=f_{|G_{S}}:G_{S}\to G$ is an
epimorphism and $\ker \varphi\subset \ker f$. Therefore, to prove
Lemma it suffices to show that $\ker f\subset G_{S}$. To show this,
let us consider the natural epimorphism $f_1:G_1\to G_2=G_1/\ker
\varphi$ and the natural epimorphism $\psi:G_2\to G$ induced by $f$.
Obviously, $\psi:(G_2,f_1(O_1))\to(G,O)$ is a cover of equipped
finite groups and $\psi_{|H}:H\to G$ is an isomorphism, where
$H=f_1(G_{S})$. Therefore $G_2\simeq \ker \psi\times G$.
Consequently, $\ker \psi =0$, since $\psi_*:H_1(G_2,\mathbb Z)\to
H_1(G,\mathbb Z)$ is
an isomorphism. \qed \\

If $S$ is a section of a cover $f:(G_1,O_1)\to (G,O)$, then Lemma
\ref{section} implies that $O_S=S^{G_1}$ is an equipment of $G_1$
and $f:(G_1,O_S)\to (G,O)$ is also a cover of equipped groups.

Below, we fix a section $S$ of a cover $f:(G_1,O_1)\to (G,O)$. Then
the cover $f$ can be considered as a quasi-cover.

In notations used in Section \ref{sec1}, consider the universal
central $C$-extension $\alpha_O:(\overline G,\overline O)\to (G,O)$
of an equipped finite group $(G,O)$. We have two natural
epimorphisms $h_O:\mathbb F_O\to G=\mathbb F_O/R_O$ and
$\beta_O:\mathbb F_O\to \overline G=\mathbb F_O/[\mathbb F_O,R_O]$
such that $h_O=\alpha_O\circ \beta_O$.

\begin{lem} \label{lem}
Let $f:(G_1,O_1)\to (G,O)$ be a quasi-cover of equipped finite
groups. Then there is an epimorphism $\alpha_{S}:(\overline
G,\overline O)\to (G_1,O_S)$ of equipped groups such that
$\alpha_O=f\circ\alpha_{S}$.
\end{lem}

\proof By Lemma \ref{section}, there is an epimorphism
$h_{S}:\mathbb F_O\to G_1$ defined by $h_{S}(x_g)=\widehat g\in S$
for all $g\in G$, where $\widehat g= f_{|S}^{-1}(g)$. Denote by
$R_{S}=\ker h_{S}$. Obviously, we have $f\circ h_{S}=h_O$. Therefore
$R_{S}\subset R_O$.

Let us show that the group $[\mathbb F_O,R_O]$ is a subgroup of
$R_{S}$. Indeed, consider any $w\in R_O$. Then, as an element of
$G_1$, the element $w\in \ker f$ and, consequently, $w$ belongs to
the center of $G_1$. In particular, it commutes with any generator
$\widehat g\in S$ of $G_1$ and hence $[w,x_g]\in R_{S}$, that is,
$[\mathbb F_O,R_O]\subset R_{S}$.

The inclusion $[\mathbb F_O,R_O]\subset R_{S}$ implies the desired
epimorphism $\alpha_{S}$. \qed \\

We say that a cover (resp., a quasi-cover) of equipped finite groups
$f:(G_1,O_1)\to (G,O)$ is {\it maximal} if for any cover of equipped
finite groups $f_1:(G_2,O_2)\to (G_1,O_1)$ such that $f_2=f\circ
f_1$ is also a cover (resp., quasi-cover) of equipped finite groups,
the epimorphism $f_1$ is an isomorphism.

\begin{thm} \label{max}
For any cover {\rm (}resp., quasi-cover{\rm)} of equipped finite
groups $f:(G_1,O_1)\to (G,O)$, there is a maximal cover {\rm
(}resp., quasi cover{\rm )} $f_2:(G_2,O_2)\to (G,O)$ for which there
is a cover $f_1:(G_2,O_2)\to (G_1,O_S)$ such that
\begin{itemize}
\item[($i$)] $f_2=f\circ f_1$;
\item[($ii$)] $\ker f_2\simeq H_2(G,\mathbb Z)$ {\rm (}resp.,
$[\overline G,\overline G]\cap \ker f_2\simeq H_2(G,\mathbb Z)${\rm
)}.
\end{itemize}
\end{thm}

\proof Consider the epimorphism $\alpha_{S}: (\overline G,\overline
O)\to (G_1,O_S)$ defined in the proof of Lemma \ref{lem}. The group
$\ker \alpha_{S}$ is a subgroup of the center of $\overline G$.

Since $(\overline G,\overline O)$ is a $C$-group and $\overline O$
consists of $M$ conjugacy classes, where $M\leq |O|=\text{rk}\,
\mathbb F_O$, then $H_1(\overline G,\mathbb Z)=\overline
G/[\overline G,\overline G]=\mathbb Z^M$. Let $\tau:\overline G\to
\mathbb Z^M$ be the natural  homomorphism (that is, $\tau$ is the
type homomorphism $\overline G\to H_1(\overline G,\mathbb Z)$, see
Introduction). The image $\tau(\ker \alpha_{S})$ is a sublattice of
maximal rank in $\mathbb Z^M$. Let us choose a $\mathbb Z$-free
basis $a_1,\dots , a_M$ in $\tau(\ker \alpha_{S})$ and choose
elements $\overline g_i\in \ker \alpha_{S}$, $1\leq i\leq M$, such
that $\tau(\overline g_i)=a_i$.

Denote by $H_{S}$ a group generated by the elements $\overline g_i$,
$1\leq i\leq M$, and denote by $K_S=[\overline G,\overline G]\cap
\ker \alpha_S$. Then it is easy to see that $H_S\simeq \mathbb Z^M$
is a subgroup of the center of $\overline G$, the intersection
$H_S\cap [\overline G,\overline G]$ is trivial, and $\ker \alpha_S
\simeq K_S\times H_S$.

Denote by $G_2=\overline G/H_S$ the quotient group and by
$\alpha_{H_S}:\overline G\to G_2$, $f_1:G_2\to G_1$ the natural
epimorphisms. We have $\alpha_S=f_1\circ \alpha_{H_S}$. Denote also
by $O_2=\alpha_{H_S}(\overline O)$. Then it is easy to see that
$\alpha_{H_S}:(\overline G,\overline O)\to (G_2,O_2)$ and
$f_1:(G_2,O_2)\to (G_1,O_S)$ are central extensions of equipped
groups.

By construction, it is easy to see that $[\overline G,\overline
G]\cap\ker \alpha_{H_S}$ is trivial and $\ker f_1\subset [G_1,G_1]$
is a subgroup of the center of $G_1$. Therefore the epimorphism
$f_1$ is a cover of equipped groups. In addition, it is easy to see
that $\alpha_O=f_1\circ\alpha_{H_S}$ and $f_2=f\circ
f_1:(G_2,O_2)\to (G,O)$ is a cover (resp., quasi-cover) of equipped
groups. We have
$$K_S\simeq\ker f_1\subset \alpha_{H_O}([\overline G,\overline G]
\cap \ker \alpha_0)=\alpha_{H_O}(H_2(G,\mathbb Z))\subset
[G_2,G_2].$$ Therefore, if $k_{f_i}=|\ker f_i|$, $i=1,2$, is the
order of the group $\ker f_i$ and $k_f$ is the order of $\ker f$,
then
\begin{equation} \label{kkk} h_2(G)=k_{f_2}=k_{f_1}k_f.\end{equation}

Since we can repeat the construction described above to the cover
(resp., quasi-cover) $f_2$ and applying again equality (\ref{kkk}),
where new $f$ is our $f_2$ and new $f_1$ is a cover existence of
which follows from assumption that old $f_2$ is not maximal, we
obtain that new $f_1$ is an isomorphism, that is, the covering $f_2$
is maximal. \qed

In the case then $f_1:(G,G\setminus
\{ 1\})\to (G,G \setminus \{1\})$ is an isomorphism of equipped
finite groups, a maximal cover $f_2: (G_2,O_2)\to (G,G\setminus \{
1\})$, constructed in the proof of Theorem \ref{max}, will be called
a {\it universal maximal cover}.

\begin{cor} \label{cor-cov} For any equipped finite group $(G,O)$
there is a maximal cover of equipped groups.

For any cover {\rm (}resp., quasi-cover{\rm )} $f:(G_1,O_1)\to
(G,O)$ of equipped finite groups, $k_f=|\ker f|\leq h_2(G)$ (resp.,
$k_f=|\ker f\cap[G_1,G_1]|\leq h_2(G)$) and $f$ is maximal if and
only if $k_f= h_2(G)$.
\end{cor}

\subsection{The ambiguity index of a quasi-cover of equipped group}
Let $(\widetilde G, \widetilde O)$ be the $C$-group associated with
an equipped group $(G,O)$ and $\beta_O:(\widetilde G, \widetilde
O)\to (G,O)$ the natural epimorphism of equipped groups (see
definitions in Section \ref{sec1}).

\begin{thm} \label{C-group} Let $f:(G_1,O_1)\to (G,O)$ be a quasi-cover
of equipped finite groups. Then there is a natural $C$-epimorphism
$\kappa_S:(\overline G,\overline O)\to (\widetilde G_1, \widetilde
O_S)$ such that $\kappa_O=\widetilde f\circ \kappa_S$ and
$\alpha_O=\beta_O\circ\widetilde f\circ \kappa_S =f\circ
\beta_{O_S}\circ \kappa_S$,  where the $C$-epimorphism
$\kappa_O:(\overline G,\overline O)\to (\widetilde G, \widetilde O)$
is defined in Section \ref{sec1} and the $C$-epimorphism $\widetilde
f: (\widetilde G_1,\widetilde O_S)\to (\widetilde G, \widetilde O)$
is associated with $f$.
\end{thm}

\proof In notations used in the proof of Lemma \ref{lem}, we have an
inclusion $R_S\subset R_O$ of normal subgroups of $\mathbb F_O$
which induces $f:G_1=\mathbb F_O/R_S\to G=\mathbb F_O/R_O$.

Let $\widetilde R_S\subset R_S$ be the normal subgroup normally
generated by the elements of $R_S$ of the form
$w_{i,j}^{-1}x_{g_i}w_{i,j}x_{g_j}^{-1}$, where $w_{i,j}\in \mathbb
F_O$ and $x_{g_i},x_{g_j}\in X_O$. For any $w\in R_O$ and any
generator $x_g$, $g\in O$, the commutator $[x_g,w]\in R_S$, since
$f$ is a central extension of groups. Therefore
\begin{equation} \label{inc}
[\mathbb F_O,R_O]\subset \widetilde R_S.\end{equation}

By Claim \ref{cl3}, $\widetilde G_1\simeq \mathbb F_S/\widetilde
R_S$. Therefore inclusion (\ref{inc}) induces an epimorphism
$\kappa_S: \overline G=\mathbb F_O/[\mathbb F_O,R_O]\to \mathbb
F/\widetilde R_S\simeq \widetilde G_1$ Obviously, the
$C$-epimorphism $\kappa_S:(\overline G,\overline O)\to (\widetilde
G_1, \widetilde O_S)$
satisfies all properties claimed in Theorem \ref{C-group}. \qed \\

Let $f:(G_1,O_1)\to (G,O)$ be a cover (resp., quasi-cover) of
equipped finite groups and $\widetilde f_S:(\widetilde
G_1,\widetilde O_S)\to (\widetilde G,\widetilde O)$ $C$-epimorphism
associated with $f:(G_1,O_S)\to (G,O)$. Denote by $k_f$ the order of
the group $\ker f\cap [G_1,G_1]$ and by $k_{\widetilde f_S}$ the
order of the group $\ker \widetilde f_S\cap [\widetilde
G_1,\widetilde G_1]$.

\begin{cor} \label{cor?} Let $f:(G_1,O_1)\to (G,O)$ be a quasi-cover
of equipped finite groups, $S$ a section of $f$. Then
$$h_2(G)=a_{(G,O)}k_{\widetilde f_S}k_S=k_fa_{(G_1,O_S)}k_S,$$
where $k_S$ is the order of the group $\ker \kappa_S\cap [\overline
G,\overline G]$.
\end{cor}

\begin{cor} \label{cor??} Let $f:(G_1,O_1)\to (G,O)$ be a cover
{\rm (}resp., quasi-cover{\rm )} of equipped finite groups, $S$ a
section of $f$. Then for any equipment $\widehat O$ of $G_1$ {\rm
(}resp., such that $O_1\subset \widehat O${\rm )} we have an
inequality $a_{(G_1,\widehat O)}\leq h_2(G)$.

If $f$ is maximal, then $a_{(G_1,\widehat O)}=1$.
\end{cor}

\proof If $f$ is a cover, then $f:(G_1,\widehat O)\to (G,f(\widehat
O))$ is also a cover of equipped groups and $a_{(G_1,\widehat
O)}\leq h_2(G)$ by Corollary \ref{cor?}.

As it was mention in the Introduction, we have $a_{(G_1,\widehat
O)}\leq a_{(G_1,O_1)}$ if $O_1\subset \widehat O$ and if $f$ is a
quasi-cover, then $a_{(G_1,O_1)}\leq h_2(G)$ by Corollary
\ref{cor?}.

If $f$ is maximal, then $k_f=h_2(G)$ by Corollary \ref{cor-cov} and
therefore if $f$ is a cover then $f:(G_1,\widehat O)\to
(G,f(\widehat O))$ is also maximal. It follows from Corollary
\ref{cor?} that $a_{(G_1,\widehat O)}=1$ in the case of maximal
covers, and $a_{(G_1,\widehat O)}\leq a_{(G_1,O_1)}=1$ in the case
of maximal quasi-covers $f$. \qed \\

Let $f:(G_1,O_1)\to (G,O)$ be a  cover of equipped
finite groups such that $f^{-1}(O)=O_1$. We say that $f$ {\it
splits} over a conjugacy class $C\subset O$ if $f^{-1}(C)$ consists
of at least two conjugacy classes of $G_1$. The number $s_{f}(C)$ of
the conjugacy classes containing in $f^{-1}(C)$ is called the {\it
splitting number} of the conjugacy class  $C$ for $f$. We say that
$f$ {\it splits completely} over $C$ if $s_f(C)=k_{f}$, where
$k_f=|\ker f|$.

Let $C$ be a conjugacy class in $G$. Consider the subgroups
$K_C\subset K_{G\setminus \{ 1\}}$ of the group
$$(R_{G\setminus \{ 1\}}\cap[\mathbb F_{G\setminus \{ 1\}},\mathbb
F_{G\setminus \{ 1\}}])/[\mathbb F_{G\setminus \{ 1\}},R_{G\setminus
\{ 1\}}]\simeq H_2(G,\mathbb Z),$$ where $K_C$ is generated by the
elements of $R_{G\setminus \{ 1\}}$ of the form $[x_{h},x_g]$, $h\in
G\setminus\{ 1\}$. Let $k_C$ be the order of the group $K_C$.

\begin{prop}\label{spl} Let $f:(G_1,O_1)\to (G,G\setminus \{ 1\})$
be a universal maximal cover of equipped finite groups and let $C$
be a conjugacy class in $G$. Then $h_2(G)=s_f(C)k_C$.\end{prop}

\proof For  $g\in C$ the preimage $f^{-1}(C)$ consists of the
conjugacy classes of the elements $zx_g$, where $z\in \ker
f=(R_{G\setminus \{ 1\}}\cap[\mathbb F_{G\setminus \{ 1\}},\mathbb
F_{G\setminus \{ 1\}}])/[\mathbb F_{G\setminus \{ 1\}},R_{G\setminus
\{ 1\}}]\simeq H_2(G,\mathbb Z)$. Note that $\ker f\subset ZG_1$ and
$\ker f$ acts transitively on the set of the conjugacy classes
$C_1,\dots C_{k_f(C)}$ involving in $f^{-1}(C)$, $z(C_i)=C_j$ if
$z\overline g\in C_j$ for $\overline g\in C_i$.

Let $x_g\in C_1$, where $g\in C$. Then $z(C_1)=C_1$ if and only if
for some $w\in G_1$ we have $wx_gw^{-1}=zx_g$, that is, $z=[w,x_g]$.

If $f(w)=h$ then $w=z_1x_{h}$ for some $z_1\in \ker f$ and therefore
$z=[x_{h},x_g]$, that is, $z\in K_C$. The inverse statement that
each element $z\in K_C$ leaves fixed the conjugacy class $C_1$ is
obvious.    \qed

\begin{prop} \label{split} Let $f:(G_1,O_1)\to (G,G\setminus \{ 1\})$ be a universal
maximal cover of equipped finite groups. Then $a_{(G,O)}=h_2(G)$ if
and only if $f$ splits completely over each conjugacy class
$C\subset O$.

If $s_f(C)=1$ for some conjugacy class $C\subset O$ then
$a_{(G,O)}=1$.
\end{prop}

\proof We have $k_f=h_2(G)$.

The map $g \mapsto x_g$ is a section in $O_1$. Denote by $\overline
O$ the equipment of $G_1$ consisting of the elements conjugated to
$x_g$, $g\in O$. Therefore $f: (G_1,\overline O)\to (G,O)$ is a
maximal cover of equipped groups and Proposition \ref{split} follows
from Corollary \ref{cor?}.  \qed

\begin{prop} \label{split1} Let $f:(G_1,O_1)\to (G,G\setminus \{ 1\})$ be a universal
maximal cover of equipped finite groups and let $C_1\subset O$ and
$C_2\subset O$ be two conjugacy classes containing in an equipment
of $G$. Then $a_{(G,O)}=1$ if $s_f(C_1)$ and $s_f(C_2)$ are coprime.
\end{prop}
\proof Follows from Corollary \ref{cor?}, since the group $\ker
\widetilde f_S\cap [\widetilde G_1, \widetilde G_1]\subset
H_2(G,\mathbb Z)$ contains two subgroups $K_{C_1}$ and $K_{C_2}$
whose indices in $H_2(G,\mathbb Z)$ are coprime. \qed

\begin{prop} \label{split2} Let $f:(G_1,O_1)\to (G,G\setminus \{ 1\})$ be a universal
maximal cover of equipped finite groups and let $h_2(G)=pq$, where
$p$ and $q$ are coprime integers. Let $C_1\subset O$ be a conjugacy
class such that $s_f(C_1)=q$ and let $s_f(C)$ is coprime with $p$
for each conjugacy class $C\subset O$. Then the ambiguity index
$a_{(G,O)}=p$.
\end{prop}
\proof Follows from Corollary \ref{cor?}, since the group $\ker
\widetilde f_S\cap [\widetilde G_1, \widetilde G_1]\subset
H_2(G,\mathbb Z)$ generated by the subgroups $K_{C_1}$ of index $p$
in $\ker f$ and subgroups of indices also coprime with $p$. \qed

\subsection{The ambiguity indices of
symmetric groups and alternating groups} \label{subsym} In
\cite{H-H}, it was proved the following theorems
\begin{thm} {\rm (Theorem 3.8 in \cite{H-H})} \label{H}
Let $\widetilde{\Sigma}_d$ be a maximal cover of the symmetric group
$\Sigma_d$. The conjugacy classes of $\Sigma_d$ which split in
$\widetilde{\Sigma}_d$ are: (a) the classes of even permutations
which can be written as a product of disjoint cycles with no cycles
of even length; and (b) the classes of odd permutations which can be
written as a product of disjoint cycles with no two cycles of the
same length {\rm (}including 1{\rm )}.
\end{thm}

\begin{thm} {\rm (Theorem 3.9 in \cite{H-H})} \label{H1}
Let $\widetilde{\mathbb A}_d$ be the maximal cover of the
alternating group $\mathbb A_d$. The conjugacy classes of $\mathbb
A_d$ which split in $\widetilde{\mathcal A}_d$ are: (a) the classes
of permutations whose decompositions into disjoint cycles have no
cycles  of even length; and (b) the classes of permutations which
can be expressed as a product of disjoint cycles with at least one
cycle of even length and with no two cycles of the same length {\rm
(}including 1{\rm )}.
\end{thm}

Remind that, by definition, an equipment $O$ of $\Sigma_d$ must
contain a conjugacy class of odd permutation since the elements of
the equipment must generate the group.

It is well known that for the symmetric group $\Sigma_d$, $d\geq 4$,
and for the alternating group $\mathbb A_d$, $d\neq 6,7$, $d\geq 4$,
the Schur multiplier $h_2(\Sigma_d)=h_2(\mathbb A_d)=2$. The
following theorems are straightforward consequences of Proposition
\ref{spl} and Theorems \ref{split} -- \ref{H1}.
\begin{thm} \label{sym} Let $O$ be an equipment of a symmetric group
$\Sigma_d$. Then $a_{(\Sigma_d,O)}=2$ if and only if $O$ consists of
conjugacy classes of odd permutations such that they can be written
as a product of disjoint cycles with no two cycles of the same
length {\rm (}including 1{\rm )} and conjugacy classes of even
permutations such that they can be written as a product of disjoint
cycles with no cycles of even length. Overwise,
$a_{(\Sigma_d,O)}=1$.
\end{thm}

\begin{thm} \label{alt} Let $O$ be an equipment of an alternating group
$\mathbb A_d$, $d\neq 6,7$. Then $a_{(\mathbb A_d,O)}=2$ if and only
if $O$ consists of  conjugacy classes of permutations whose
decompositions into disjoint cycles have no cycles  of even length
and the classes of permutations which can be expressed as a product
of disjoint cycles with at least one cycle of even length and with
no two cycles of the same length {\rm (}including 1{\rm )}.
Overwise, $a_{(\mathbb A_d,O)}=1$.
\end{thm}

It is well known that in the case when $d= 6,7$, the Schur
multiplier $h_2(\mathbb A_d)=6$.

For $\sigma\in \mathbb A_d$ denote by $c(\sigma)=(l_1,\dots, l_m)$
the cycle type of permutation $\sigma$, that is, the collection of
lengths $l_i$ of non-trivial (that is $l_i\geq 2$) cycles entering
into the factorization of $\sigma$ as a product of disjoint cycles.
For a conjugacy class $C$ in $\mathbb A_d$ the collection
$c(C)=c(\sigma)$ is called the {\it cycle type} of $C$ if $\sigma\in
C$. It is well known that the cycle type $c(C)$ does not depend on
the choice of $\sigma\in C$ and there are at most two conjugacy
classes in $\mathbb A_d$ of a given cycle type $c$.

The group $\mathbb A_d$, $d=6,7$, has the following non-trivial
conjugacy classes:
\begin{itemize}
\item[(I)] two  conjugacy classes of each cycle type  $(5)$,
$(2,4)$, and (if $d=7$) $(7)$;
\item[(II)] two conjugacy classes of cycle type $(3)$ and one
conjugacy class of cycle type $(3,3)$;
\item[(III)] one conjugacy class of cycle type $(2,2)$ and one
conjugacy class of cycle type $(2,2,3)$ if $d=7$.
\end{itemize}
\begin{prop} \label{d67} The ambiguity index $a_{(\mathbb A_d,O)}$, $d=6,7$,
takes the following values:
\begin{itemize}
\item[(I)] $a_{(\mathbb A_d,O)}=6$ if $O$ contains only the elements
of conjugacy classes of type {\rm (I)};
\item[(II)] $a_{(\mathbb A_d,O)}=2$ if $O$ contains only the elements
of conjugacy classes of type {\rm (I)} and the elements of at least
one conjugacy class of type {\rm (II)};
\item[(III)] $a_{(\mathbb A_d,O)}=3$ if $O$ contains only the elements
of conjugacy classes of type {\rm (I)} and the elements of at least
one conjugacy class of type {\rm (III)};
\item[(II+III)] $a_{(\mathbb A_d,O)}=1$ if $O$ contains the elements
of at least one  conjugacy class of type {\rm (II)} and the elements
of at least one conjugacy class of type {\rm (III)}.
\end{itemize}
\end{prop}

\proof Let $f:(G_1,O_1)\to(\mathbb A_d,\mathbb A_d\setminus\{ 1\})$
be the universal maximal cover.

Note that, by \cite{Kun1}, $a_{((\mathbb A_d,\mathbb A_d\setminus\{
1\})}=1$. Therefore  there exist elements
$\sigma_1,\dots,\sigma_4$ in $\mathbb A_d$ such that 
$[x_{\sigma_1},x_{\sigma_2}]$ and $[x_{\sigma_3},x_{\sigma_4}]$  in
$([\mathbb F_{\mathbb A_d\setminus\{ 1\}},\mathbb F_{\mathbb
A_d\setminus\{ 1\}}]\cap R_{\mathbb A_d})/[\mathbb F_{\mathbb
A_d\setminus\{ 1\}},R_{\mathbb A_d}]$ have, respectively, order two
and three.

It is easy to see that for an element $\sigma$ belonging to a
conjugacy class $C$ of type (I) the centralizer $Z(\sigma)\subset
\mathbb A_d$ of the element $\sigma$ is a cyclic group generated by
$\sigma$. Therefore $K_C$ is the trivial group and hence
$s_f(C)=h_2(\mathbb A_d)$. Therefore, by Proposition \ref{split},
$a_{(\mathbb A_d,O)}=6$ if $O$ contains only the elements of
conjugacy classes of type {\rm (I)}.

Let $\sigma$ is of cycle type $(2,2,3)$. Without loss of generality,
we can assume that $\sigma=\sigma_1\sigma_2$, where
$\sigma_1=(1,2)(3,4)$ and $\sigma_2=(5,6,7)$. Then the centralizer
$Z(\sigma)\subset \mathbb A_d$ of $\sigma$ is $\text{Kl}_4\times
\langle \sigma_2\rangle$, where
$\text{Kl}_4=\langle\sigma_1\rangle\times \langle\sigma_3\rangle$
and $\sigma_3=(1,3)(2,4)$. We have
$[x_{\sigma},x_{\sigma_3,\sigma_2^{\pm
1}}]=[x_{\sigma_1},x_{\sigma_3}]$ in the group $\mathbb F_{\mathbb
A_d\setminus\{ 1\}}/[\mathbb F_{\mathbb A_d\setminus\{
1\}},R_{\mathbb A_d}]$. Therefore $K_C$, where $C$ has type
$(2,2,3)$, is a group of order at most two since the order of
$\sigma_1$ is two (see Lemma \ref{order}) and it is of order two if
and only if $[x_{\sigma_1},x_{\sigma_3}]$ is not the unity in
$\mathbb F_{\mathbb A_d\setminus\{ 1\}}/[\mathbb F_{\mathbb
A_d\setminus\{ 1\}},R_{\mathbb A_d}]$. But, the embeddings
$\langle\sigma_1,\sigma_3\rangle\subset \mathbb A_d\subset\Sigma_d$
define a sequence  of homomorphisms
$$H_2(\langle\sigma_1,\sigma_3\rangle,\mathbb Z)\to H_2(\mathbb
A_d,\mathbb Z)\to H_2(\Sigma_d,\mathbb Z)$$ such that the image of
the non-trivial element $[x_{\sigma_1},x_{\sigma_3}]$ in
$H_2(\langle\sigma_1,\sigma_3\rangle,\mathbb Z)$ is a non-trivial in
$H_2(\Sigma_d,\mathbb Z)$. Therefore $s_f(C)=3$ for the conjugacy
class $C$ of cyclic type $(2,2,3)$ and, similarly, $s_f(C)=3$ for
the conjugacy class $C$ of cyclic type $(2,2)$, since $K_C$ is a
subgroup of $H_2(\mathbb A_d,\mathbb Z)\simeq \mathbb Z/6\mathbb Z$
generated by the elements of the second order (see Proposition
\ref{spl}) and only the elements of $K_{C_1}$ and $K_{C_2}$ can
generate the subgroup of order two in $H_2(\mathbb A_d,\mathbb Z)$.

Let $\sigma$ is of cycle type $(3,3)$. Without loss of generality,
we can assume that $\sigma=\sigma_1\sigma_2$, where
$\sigma_1=(1,2,3)$ and $\sigma_2=(4,5,6)$. Then the centralizer
$Z(\sigma)\subset \mathbb A_d$ of $\sigma$ is
$\langle\sigma_1\rangle\times \langle \sigma_2\rangle$. Therefore
$[x_{\overline \sigma},x_{\sigma}]$ is not the unity in $\mathbb
F_{\mathbb A_d\setminus\{ 1\}}/[\mathbb F_{\mathbb A_d\setminus\{
1\}},R_{\mathbb A_d}]$ only if  $\overline\sigma=\sigma_1^{\pm 1}$,
either $\overline\sigma=\sigma_2^{\pm 1}$, or $\overline \sigma=
\sigma_1\sigma_2^{-1}$, or $\overline \sigma=
\sigma_1^{-1}\sigma_2$. We have
$$[x_{\sigma_1\sigma_2^{-1}},x_{\sigma}]=[x_{\sigma_1},x_{\sigma_2}]
[x_{\sigma_2^{-1}},x_{\sigma_1}]=[x_{\sigma_1},x_{\sigma_2}]^2$$ in
$\mathbb F_{\mathbb A_d\setminus\{ 1\}}/[\mathbb F_{\mathbb
A_d\setminus\{ 1\}},R_{\mathbb A_d}]$ and, similarly,
$[x_{\sigma^{-1}\sigma_2},x_{\sigma}]=[x_{\sigma_1},x_{\sigma_2}]$,
since the elements $x_{\sigma}x_{\sigma_2}^{-1}x_{\sigma_1}^{-1}$,
$x_{\sigma_1\sigma_2^{-1}}x_{\sigma_2}x_{\sigma_1}^{-1}$ belong to
the center of the group $\mathbb F_{\mathbb A_d\setminus\{
1\}}/[\mathbb F_{\mathbb A_d\setminus\{ 1\}},R_{\mathbb A_d}]$.
Therefore the group $K_{C_1}$ is a nontrivial group of order three
if and only if $K_{C_2}$ is a nontrivial group of order three, where
$C_1$ is a conjugacy class of the cycle type $(3)$ and $C_2$ is the
conjugacy class of the cycle type $(3,3)$, and hence
$s_f(C_1)=s_f(C_2)=2$. Now Proposition \ref{d67} follows from
Propositions \ref{split} -- \ref{split2}. \qed

\section{Cohomologycal description of the ambiguity indices} \label{sec4}
In notations used in Section \ref{sec1}, for an equipped finite
group $(G,O)$ a subgroup $K_{(G,O)}$ of $H_2(G,\mathbb Z)$ was
defined as follows: $K_{(G,O)}$ is the subgroup of $(R_O\cap[\mathbb
F_O,\mathbb F_O])/[\mathbb F_O,R_O]$ generated by the elements of
$R_O$ of the form $[w,x_g]$, where $g\in O$ and $w\in \mathbb F_O$,
and $k_{(G,O)}$ is it's order.

Denote $$B_{(G,O)}=K_{(G,O)}^{\bot}=\{ \varphi \in Hom
(H_2(G,\mathbb Z),\mathbb Q/\mathbb Z)\mid \varphi(w)=0\,\,
\text{for all}\, \, w\in K_{(G,O)}\}$$ a subgroup of $H^2(G,\mathbb
Q/\mathbb Z)$ dual to $K_{(G,O)}$. As in the proof of Theorem
\ref{main}, it is easy to show that
$$B_{(G,O)})=\ker[H^2(G,\mathbb Q/\mathbb
Z)\to \bigotimes_{A\subset G}H^2(A,\mathbb Q/\mathbb Z)],$$ where
$A$ runs over all abelian subgroups of $G$ generated by two elements
$g\in O$ and $h\in G$. Let $b_{(G,O)}$ be the order of the group
$B_{(G,O)}$. In particular, $b_{(G,G\setminus \{ 1\})}=b_0(G)$.

The next theorem immediately  follows from Theorem \ref{thm2}.
\begin{thm} \label{thm?} For an equipped finite group $(G,O)$
we have $a_{(G,O)}=b_{(G,O)}$.
\end{thm}

The group $H^2(G,Q/Z)$ is a  direct sum of primary components
$H^2(G,Q/Z) = \Sigma_p H^2(G,Q/Z)_{p}$ where primes $p$ run through
a subset of primes dividing the order of of $H^2(G,Q/Z)$ and hence
$G$. Therefore we have the following:

\begin{prop} \label{prop3} If the set of conjugacy classes
$O$ consists of all classes of power of prime order then
$a_{(G,O)}=b_0(G)$. Moreover it is sufficient to consider such
classes only for primes dividing $h_2(G)$.
\end{prop}

Note that $H^2(G,\mathbb Q/\mathbb Z)_{p}$ embeds into
$H^2(Syl_p(G),\mathbb Q/\mathbb Z)_{p}$ where $Syl_p(G)$ is a Sylow
$p$-subgroup of $G$. Similarly the $p$-primary component
$B_0(G)_{p}$ is a subgroup of $B_0(Syl_p(G))$.

More explicit versions of Proposition \ref{prop3} for different
groups provide with simple methods to compute $B_O(G)$

\section{An example of a finite group $G$ with $b_0(G)>1$}
\label{sec5} The following groups where constructed in the article
of Saltman \cite{Sa84}.

Consider a finite $p$-group $G_p$ which is a central extension of
$\mathbb Z_p^4= A_p$ with generators $x_i$. The center of $G_p$ is
generated by pairwise commutators
 $x_i x_j x_i^{-1} x_j^{-1}= [x_i,x_j]$ with one relation
between $[x_1,x_2] [x_3,x_4] =1$. Thus there is natural exact
sequence:
$$ 1\to \mathbb Z_p^5 \to G_p \to  A_p \to 1 $$.
\begin{lem} {\rm (}\cite{BO87}, \cite{Sa84}{\rm )}  $B_o(G_p)= \mathbb Z/p$.
\end{lem}
\proof It is shown in \cite{BO87} using standard spectral sequence
that for a central extension $G$ of an abelian group $A$ the group
$B_0(G)$ is contained in the image of $H^2(A,\mathbb Q/\mathbb Z)$
in $H^2(G,\mathbb Q/\mathbb Z)$. The group $H^2(A_p,\mathbb
Q/\mathbb Z) = \mathbb Z_p^6$ which is generated by elements
$[x_i,x_j]^*$. The kernel of the map $H^2(A_p,\mathbb Q/\mathbb
Z)\to H^2(G_p,\mathbb Q/\mathbb Z)$ is naturally dual to the center
$\mathbb Z_p^5$ of $G_p$. Thus the image of $H^2(A_p,\mathbb
Q/\mathbb Z)$ in  $H^2(G_p,\mathbb Q/\mathbb Z)$ is a cyclic
$p$-group generated by one element $w$. Let us show that the latter
is in $B_0(G_p)$. It is enough to check that it is trivial on any
abelian subgroup in $G_p$ which surjects onto rank $2$ subgroup
$\mathbb Z_p^2\subset \mathbb Z_p^4= A_p$. However $G_p$ does not
contain such subgroups. Indeed assume that the restriction $w$ on a
subgroup with generators $x_1,y_1\in G_p$ is trivial. It means that
the  commutator $[x,y] =1$ in $G_p$ where $x,y$ are projections of
$x_1,y_1$ into $A_p$. On the other hand the only nontrivial relation
between commutators of elements in $A_p$ is  $[x_1,x_2] [x_3,x_4]
=1$ which is not equal to $[x,y]$ for any pair $x,y\in A_p$. Hence
$w$ restricts trivially onto any subgroup with two generators in
$G_p$ and generates $B_0(G_p)$.

 \ifx\undefined\bysame
\newcommand{\bysame}{\leavevmode\hbox to3em{\hrulefill}\,}
\fi

\end{document}